# A SPATIALLY EXPLICIT MODEL FOR COMPETITION AMONG SPECIALISTS AND GENERALISTS IN A HETEROGENEOUS ENVIRONMENT

BY N. LANCHIER AND C. NEUHAUSER

*University of Minnesota*

Competition is a major force in structuring ecological communities. The strength of competition can be measured using the concept of a *niche*. A niche comprises the set of requirements of an organism in terms of habitat, environment and functional role. The more niches overlap, the stronger competition is. The niche breadth is a measure of specialization: the smaller the niche space of an organism, the more specialized the organism is. It follows that, everything else being equal, generalists tend to be more competitive than specialists. In this paper, we compare the outcome of competition among generalists and specialists in a spatial versus a nonspatial habitat in a heterogeneous environment. Generalists can utilize the entire habitat, whereas specialists are restricted to their preferred habitat type. We find that although competitiveness decreases with specialization, specialists are more competitive in a spatial than in a nonspatial habitat as patchiness increases.

**1. Introduction.** An ecological community is an assembly of populations that inhabit the same spatial location at the same time. Understanding the forces that shape such communities has been the goal of community ecology.

Competition was identified early on as an important driver of community assembly and is closely linked to the concept of a *niche*, which was introduced by Elton in 1927 [12]. Elton defined niche as "the status of an animal in its community" and used it to explain how multiple species can coexist within a community. Hutchinson [15] formalized this concept and defined niche mathematically as a subset of an $n$-dimensional hypervolume. Species with overlapping niches compete with each other. One of the main insights has

Received August 2005; revised December 2005.
Supported in part by NSF DMS-00-83498.

*AMS 2000 subject classifications.* Primary 60K35; secondary 82C22.

*Key words and phrases.* Multitype contact process in heterogeneous environment, specialist, generalist.

been that in order for species to coexist, they cannot occupy identical niches, that is, identical subsets of the niche space. Although delineating the exact niche of a species is nearly impossible in nature, the concept has proven quite useful, particularly in mathematical models where niches can be identified more easily.

Species that coexist are able to divide up the niche space so that each species occupies some subset of the niche space by itself. It thus follows that communities composed of species with narrower niches can be more species-rich than those composed of species with broader niches (see, e.g., [6]). We call the former *specialists* and the latter *generalists*.

Much research has been done to understand the roles of generalists and specialists in ecological communities. No clear patterns have emerged. Some taxa consist of mostly specialists. For instance, about 90% of herbivorous insects are highly specialized, feeding on three or fewer plant families [2]. Others are dominated by generalists. For instance, 60% of human pathogens are zoonotic, that is, they infect both humans and nonhuman animals, and often infect hosts from different orders or classes [22].

Evolutionary ecologists have used optimal behaviors of individuals to gain insights into the question of relative success of specialists versus generalists, in order to assess these strategies in an evolutionary context and ultimately to predict community structure. Extending the work of Rosenzweig [21], Brown [5] found in a theoretical study on habitat selection that in an environment with two habitat types, three scenarios are possible: (a) two specialists, one for each habitat type, (b) one generalist who occupies both habitat types and (c) one generalist and one specialist. Scenario (a) occurs when there is no cost to habitat selection, whereas scenario (b) occurs when there is a cost. Scenario (c) occurs when there is an asymmetry in the two habitat types, either in terms of relative abundances or productivities.

The theoretical work mentioned above relied on nonspatial models and thus did not incorporate a way to account for the often very localized interactions between competitors. Past research on competition has demonstrated that the inclusion of a spatial component with local interactions can alter predictions based on nonspatial models [19] or offer new insights into consequences of local interactions [1]. In this paper we employ the multitype contact process [18] in a heterogeneous environment to investigate how local, competitive interactions among specialists and generalists affect their abilities to coexist. We define the environment as a patchwork of two host types on which specialists and generalists feed. Specialists feed on one of the two hosts, whereas the generalists can consume either host. We assume that the environment does not change over time.

We begin by defining the environment. To set the spatial configuration of the host population, we let $L \geq 1$ be an integer, which we will refer to as the *space scale*, introduce the $d$-dimensional box $\mathcal{H}_L = [-L, L)^d$ and define $\mathcal{H}_1$ and $\mathcal{H}_2$ as a partition of $\mathbb{Z}^d$ given by

$$\mathcal{H}_1 = \bigcup_{z \in H_1} (2Lz + \mathcal{H}_L) \quad \text{and} \quad \mathcal{H}_2 = \bigcup_{z \in H_2} (2Lz + \mathcal{H}_L),$$

where $H_1$ and $H_2$ denote the sets of $z \in \mathbb{Z}^d$ such that $z_1 + \cdots + z_d$ is even and odd, respectively. In our biological context, $\mathcal{H}_i$ will represent the set of sites that are permanently occupied by a host of type $i$, $i = 1, 2$.

Into this host population we introduce two specialist species and one generalist species that differ in their abilities to reproduce. We will refer to these species collectively as *consumers*. To describe their evolution, we consider a continuous-time Markov process whose state at time $t$ is a function $\xi_t : \mathbb{Z}^d \longrightarrow \{0, 1, 2, 3\}$. A site $x \in \mathbb{Z}^d$ is said to be occupied by a specialist of type 1 (resp. 2) if $\xi(x) = 1$ (resp. 2), and a generalist if $\xi(x) = 3$, with the condition $\xi(x) = 0$ denoting the absence of a consumer. The local dynamics at site $x$ are described by the following transition rates:

$$0 \to 1 \quad \text{at rate } \alpha \sum_{0 < \|x-z\| \leq R} \mathbb{1}_{\{\xi(z)=1; x \in \mathcal{H}_1\}} \qquad 1 \to 0 \quad \text{at rate } 1,$$

$$0 \to 2 \quad \text{at rate } \alpha \sum_{0 < \|x-z\| \leq R} \mathbb{1}_{\{\xi(z)=2; x \in \mathcal{H}_2\}} \qquad 2 \to 0 \quad \text{at rate } 1,$$

$$0 \to 3 \quad \text{at rate } \beta \sum_{0 < \|x-z\| \leq R} \mathbb{1}_{\{\xi(z)=3\}} \qquad 3 \to 0 \quad \text{at rate } 1,$$

where $\|x\| = \sup_{i=1,2,\ldots,d} |x_i|$. We refer to $R$ as the *dispersal range*. To understand the evolution of our stochastic process when $d = 2$, let us imagine an infinite chessboard where each of its squares has length $2L$ (and thus $2L \times 2L$ sites). The black and white squares represent parts of the habitat that are permanently occupied by hosts of type 1 and 2, respectively. Particles of type 3, that is, the generalists, do not see the colors of the chessboard and, in the absence of particles of type 1 and 2, perform a (translation-invariant) contact process. That is, each 3 tries to give birth onto each of its neighboring sites at rate $\beta$, the birth occurring if and only if the offspring is sent to a site that is not already occupied by a consumer, regardless of the type of the host present at that site. The dynamics of the specialists are quite different since black (resp. white) squares present unsuitable habitat for particles of type 2 (resp. 1). Specifically, if the offspring of a 1 (resp. 2) is sent to a white (resp. black) square, then the birth is suppressed.

*The mean-field model.* Before we describe the behavior of the spatially explicit stochastic model, we will look at the mean-field model [9], that is, we will pretend that all sites are independent and that the probability of a site being occupied by a consumer of type $i$ depends only on the host type at that site, but is otherwise spatially homogeneous. This then results in a system of differential equations for the densities of specialists and generalists. Let $v_{ij}$ denote the density of hosts of type $i$ associated with a consumer of type $j$, where $i = 1, 2$ and $j = 1, 2, 3$. We let

$$u_1 = \tfrac{1}{2} - v_{11} - v_{13} \quad \text{and} \quad u_2 = \tfrac{1}{2} - v_{22} - v_{23}$$

denote the densities of unassociated hosts of type 1 and 2, respectively. The constant $1/2$ corresponds to the limit of the proportions of the two host types in a spatial box whose size tends to infinity. The mean-field limit is obtained by letting the dispersal range $R$ tend to infinity. To obtain a meaningful limit, we need to rescale the parameters $\alpha$ and $\beta$ by the neighborhood size, that is, we set

$$\alpha = \frac{a}{\nu_R} \quad \text{and} \quad \beta = \frac{b}{\nu_R},$$

where $\nu_R$ denotes the size of the dispersal neighborhood, that is, $\nu_R = |\{z : 0 < \|x - z\| \leq R\}|$. The following system of differential equations describes the mean-field behavior:

$$\frac{dv_{11}}{dt} = -v_{11} + a u_1 v_{11},$$

$$\frac{dv_{22}}{dt} = -v_{22} + a u_2 v_{22},$$

$$\frac{dv_{13}}{dt} = -v_{13} + b u_1 (v_{13} + v_{23}),$$

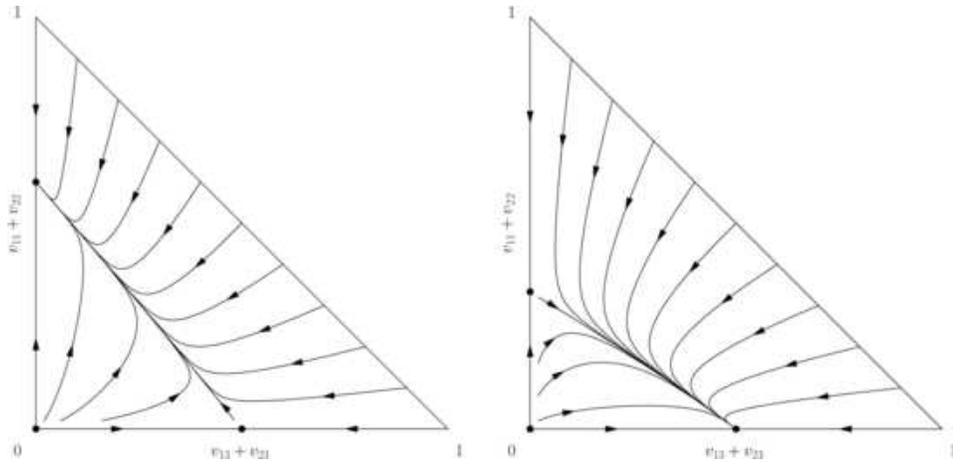

FIG. 1. *Solution curves of the mean-field model. Picture on the left: $a = 5$ and $b = 2$. Picture on the right: $a = 3$ and $b = 2$. In each of the pictures, we start with $v_{11} = v_{22}$ and $v_{13} = v_{23}$.*

$$\frac{dv_{23}}{dt} = -v_{23} + bu_2(v_{13} + v_{23}).$$

A standard linear stability analysis reveals the equilibrium behavior [20]. If $b < 1$ and $a < 2$, then the trivial equilibrium ($v_{11} = v_{22} = v_{13} = v_{23} = 0$) is the only locally stable equilibrium. If $a > 2b$ and $a > 2$, then there is a locally stable equilibrium point with $v_{11} = v_{22} > 0$ and $v_{13} = v_{23} = 0$ (Figure 1; picture on the left). If $a < 2b$ and $b > 1$, then there is a locally stable equilibrium in which $v_{11} = v_{22} = 0$ and $v_{13} = v_{23} > 0$ (Figure 1; picture on the right). Furthermore, if one of the specialists is initially absent (say, $v_{22} = 0$ at time 0), then there is an additional equilibrium in which $v_{11}$, $v_{13}$ and $v_{23}$ are positive, provided $a > 2$ and $\frac{2a}{a+2} < b < \frac{a}{2}$. This equilibrium can, however, be invaded by specialist 2, which then outcompetes the generalist, resulting in the equilibrium in which $v_{11} = v_{22} > 0$ and $v_{13} = v_{23} = 0$. There is no locally stable equilibrium in which both specialists and the generalist coexist (i.e., all three have positive densities). A necessary condition for such an equilibrium is $a = 2b > 2$, but the resulting equilibrium is not locally stable (it is neutrally stable).

*The spatially explicit model.* We now return to the spatially explicit model with parameters $\alpha$ and $\beta$. Our results show that the inclusion of a spatial structure in the form of local interactions can significantly modify the outcome of competition among specialists and generalists. While in the mean-field model, both the dispersal range $R$ and the space scale $L$ lose their meaning, the study of the spatial model reveals that setting $L$ much larger than $R$ helps the specialists to compete.

As explained above, in the absence of specialists, the process for the generalists reduces to a translation-invariant contact process with parameter $\beta$. In particular, there exists a critical value $\lambda_c \in (0, \infty)$ such that the following holds: If $\beta \leq \lambda_c$, the process converges in distribution to the "all 0" configuration. If $\beta > \lambda_c$, there exists a stationary measure $\mu_3$ that concentrates on configurations with infinitely many 3's ([3] or [17], Theorem 2.25). A simple coupling argument shows that specialists die out when $\alpha \leq \lambda_c$. To show this, we consider an initial configuration $\xi_0$ with $\xi_0(x) \neq 3$ for all $x \in \mathbb{Z}^d$. Let $\zeta_t$ denote the contact process with parameter $\alpha$ and initial configuration $\zeta_0(x) = 0$ if $\xi_0(x) = 0$, and $\zeta_0(x) = 3$ otherwise, that is, $\zeta_0$ can be deduced from $\xi_0$ by replacing each specialist by a generalist. Both processes can then be constructed graphically from the same collection of independent Poisson processes [13] in such a way that $\zeta_t$ has more consumers than $\xi_t$ at any time $t \geq 0$ (see Section 2 for a detailed description of this graphical representation). Since $\zeta_t$ (the process for the generalists) dominates $\xi_t$ (the process for the specialists) and $\zeta_t$ converges to the "all 0" configuration for $\alpha \leq \lambda_c$, it follows that $\xi_t$ also converges to the "all 0" configuration for $\alpha \leq \lambda_c$. Hence specialists die out when $\alpha \leq \lambda_c$. It is more difficult to show that specialists can survive when $\alpha > \lambda_c$ (see Theorem 3).

To investigate the competition between specialists and generalists, we now consider the process starting with infinitely many consumers of each type. Relying on a new coupling argument, we may run our stochastic process $\xi_t$ and the 3-color multitype contact process $\eta_t$ (in which particles of type 1 and 2 give birth at rate $\alpha$ and particles of type 3 give birth at rate $\beta$) on the same probability space in such a way that, starting from the same initial configuration, $\xi_t$ has more 3's and fewer 1's and 2's than $\eta_t$. According to Theorem 1 in [18] the process $\eta_t$ converges in distribution to $\mu_3$, provided $\beta > \alpha > \lambda_c$. The above-mentioned coupling then implies that the generalists outcompete the specialists. An analysis of the dual process will show that the result still holds when $\alpha = \beta > \lambda_c$. That is, while coexistence occurs for the multitype contact process in a homogeneous environment when $\alpha = \beta$ and $d \geq 3$ (see [18]), the process in this heterogeneous environment gives the generalists an advantage over specialists in any dimension. This is summarized in the following theorem where "$\Rightarrow$" denotes weak convergence, and $\mu_3$ is the upper invariant measure of the basic contact process introduced above.

THEOREM 1.  *Assume that $L$ and $R$ are fixed and that $\beta \geq \alpha$ and $\beta > \lambda_c$. If $\xi_0$ is translation-invariant and $P(\xi_0(x) = 3) > 0$, then $\xi_t \Rightarrow \mu_3$.*

This result corresponds to the case $b \geq a$ and $b > 1$ in the mean-field model and tells us that the fragmentation of the environment allows the generalists to outcompete the specialists in both the spatial and the nonspatial models, provided the birth rate of the generalists is at least as high as that of the specialists.

The next step is to investigate the behavior of the spatially explicit model when the specialists have a higher birth rate than the generalists. In the mean-field model, we have seen that the condition $a \geq 2b$ is required for the specialists to invade the generalists in their equilibrium. The behavior of the stochastic process, however, is different. Let us choose a space scale $L$ much larger than the range of interactions $R$ so that well inside black squares, specialists of type 1 and generalists behave nearly like a multitype contact process. In particular, specialists locally outcompete generalists, provided $\alpha > \beta$. The survival of specialists, however, is not clear since the contact process on a finite set (here the $2L \times 2L$ square when $d = 2$) always converges to its absorbing state, the "all 0" configuration. To show that the specialists actually survive on the infinite chessboard (with large squares), we will prove that, with probability close to 1, the 1's can invade the four diagonally adjacent black squares before going extinct. The trickiest point of the proof will be to show that, with high probability, a specialist of type 1 in the middle of a black square will produce four "invasion paths" to bring its offspring within a finite distance (independent of $L$) of each of the four corners of the square. Since white squares are forbidden for the 1's, the invasion paths will have to be contained in the black square. Comparing the set of black squares occupied by (a significant quantity of) 1's with the set of wet sites in a one-dependent oriented percolation process will then yield the following:

THEOREM 2. *Assume that $R$ is fixed and that initially there are infinitely many 1's and 3's, but no 2's. If $\alpha > \beta > \lambda_c$ and $d \geq 2$, then 1's and 3's coexist, provided $L$ is sufficiently large.*

Since specialists of type 1 and 2 do not interfere with each other, the proof of Theorem 2 also implies the following:

THEOREM 3. *Assume that $R$ is fixed and $d \geq 2$. If $\alpha > \lambda_c$ and $\alpha > \beta$, and initially there are infinitely many 1's and 2's, then 1's and 2's coexist, provided the space scale $L$ is sufficiently large.*

Unfortunately, Theorem 3 does not tell us if all three consumers coexist when $d \geq 2$. Numerical simulations on a $200 \times 200$ torus indicate that the generalists persist for a very long time along the boundaries of the squares of the chessboard, namely where the density of specialists is low, but proving or disproving that they can percolate on the infinite lattice would be difficult.

*Comparison of the spatially explicit and the mean-field models.* In the mean-field model, that is, in the absence of spatial structure, the condition for the specialists to invade the generalists is given by $a \geq 2b$. The factor 2 comes from the fact that only half of the habitat is available for each specialist. In the spatially explicit model, the dimension $d$ as well as the space scale $L$ and the dispersal range $R$ can profoundly modify the outcome of the competition. In the one-dimensional case, the behavior of the spatial model is quite easy to predict when $R \leq 2L$. Since all the offspring of a specialist are confined to their parents' segment of length $2L < \infty$ and are unable to reach other segments (because $R \leq 2L$), the specialists eventually go extinct so that the only nontrivial stationary distribution of the process is $\mu_3$. To describe the behavior of the process when $d \geq 2$, let $\alpha_c$ denote the infimum of $\alpha$'s such that the specialists can survive on the chessboard. Theorems 1 and 3 tell us that if the birth rate of the generalists is equal to $\beta$, then $\beta < \alpha_c < \infty$. Moreover, $\alpha_c \to \beta$ as the space scale $L$ tends to infinity. That is, taking large $L$ makes the specialists much more competitive in the spatially explicit model than in the mean-field model. The right-hand picture in Figure 2 illustrates the effects of $L$ on the outcome of the competition.

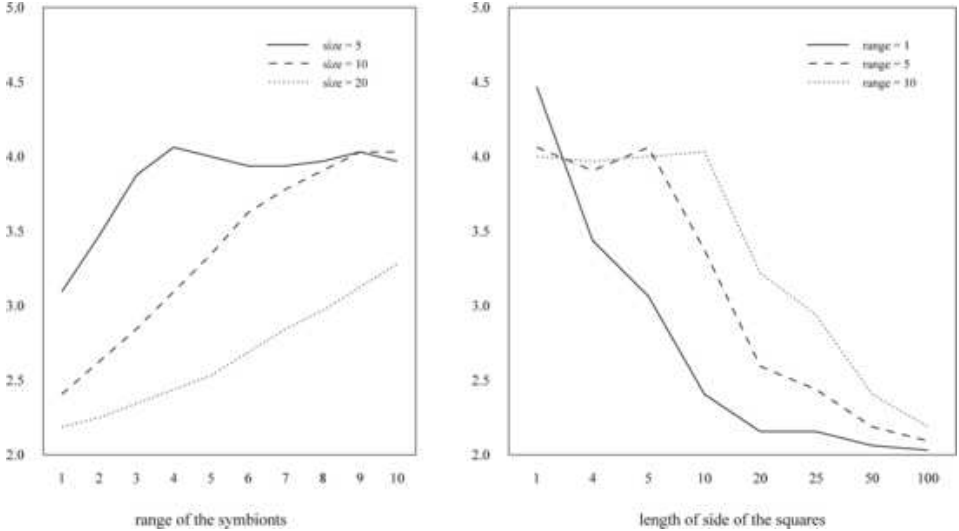

FIG. 2. *Simulation results for the spatially explicit model with parameter $\beta = 2$ on the $200 \times 200$ square with periodic boundary conditions. In both pictures, the vertical axis denotes the common birth rate of specialists. The curves represent the smallest $\alpha$ such that the density of specialists is greater than the density of generalists after 100 units of time. The simulation results suggest that when the range of the interactions of the symbionts is larger than the length of the squares, the competitiveness of the specialists reduces significantly, just as in the mean-field model. The competitiveness of the specialists increases when the range interactions is short and the habitat is coarse-grained as shown in Theorem 2.*

Another interesting corollary of Theorem 3 is that, in the absence of generalists, $\lambda_c < \alpha_c < \infty$ and $\alpha_c \to \lambda_c$ as $L$ tends to infinity. In conclusion, our (multitype) contact process on the chessboard behaves nearly like the corresponding (multitype) contact process in a homogeneous environment when the size of the squares is large, while their respective nonspatial versions differ significantly.

**2. Duality. Proof of Theorem 1.** We saw in the Introduction that the behavior of the process when $\alpha < \beta$ is a straightforward consequence of Theorem 2 of [18], so we will now focus on the case $\alpha \geq \beta$. The key to the proof of Theorem 1 is duality. To define the dual process of our heterogeneous multitype contact process, we start by introducing a graphical representation from which the process with parameters $\alpha$ and $\beta$, where $\alpha \geq \beta$, can be constructed.

*The graphical representation when $\alpha \geq \beta$.* For any $x, z \in \mathbb{Z}^d$, $0 < \|x - z\| \leq R$, let $\{T_n^{x,z} : n \geq 1\}$ and $\{U_n^x : n \geq 1\}$ denote the arrival times of independent Poisson processes with parameters $\alpha$ and 1, respectively. At times $T_n^{x,z}$, we draw an arrow from site $x$ to site $z$, toss a coin with success probability $(\alpha - \beta)/\alpha$ and, if there is a success, label the arrow with an **s** (for specialist). Moreover, if $x \in \mathcal{H}_i$ and $z \in \mathcal{H}_j$ with $i \neq j$, then we label the arrow with a **g** (for generalist). This can be interpreted as follows. If, at time $T_n^{x,z}$, the site $x$ is occupied by a consumer and $z$ is empty, then the consumer at site $x$, regardless of type, gives birth to a consumer of the same type through this arrow, unless the arrow is labeled with an **s** or a **g**. Only specialists (resp. generalists) can give birth if the label is **s** (resp. **g**). We stipulate that arrows labeled with both **s** and **g** cannot be used by either consumer. To complete the construction, we put a "$\times$" at $(x, U_n^x)$ to indicate that a consumer, regardless of its type, dies.

*The dual process.* The dual process in the heterogeneous environment is defined as for the process in a homogeneous environment. We say that there is a *path* from $(x, s)$ to $(z, T)$, where $0 \leq s \leq T$, if there is a sequence of times $s_0 = s < s_1 < \cdots < s_{n+1} = T$ and a sequence of spatial locations $x_0 = x, x_1, \ldots, x_n = z$ such that the following two conditions hold:

1. for $i = 1, 2, \ldots, n$, there is an arrow from $x_{i-1}$ to $x_i$ at time $s_i$, and
2. for $i = 0, 1, \ldots, n$, the vertical segments $\{x_i\} \times (s_i, s_{i+1})$ do not contain any $\times$'s.

A *dual path* from $(x, T)$ to $(z, T - s)$ then refers to a path from $(z, T - s)$ to $(x, T)$. The *dual process starting at* $(x, T)$ is defined by setting

$$\hat{\xi}_s^{(x,T)} = \{z \in \mathbb{Z}^d : \text{there is a dual path from } (x, T) \text{ to } (z, T - s)\}.$$

The dual process is naturally defined only for $0 \leq s \leq T$ but it is convenient to assume that the Poisson processes involved in the graphical representation are defined for negative times so that the dual process can be defined for all $s \geq 0$.

We introduce the dual process in order to deduce the state of site $x$ at time $T$ from the configuration at earlier times. The aim is to understand the ancestry of $(x, T)$ by working backwards in time. As for the multitype contact process, at any time $s \geq 0$, the elements of $\hat{\xi}_s^{(x,T)}$, called *ancestors*, can be arranged according to the order in which they determine the type of $(x, T)$. This hierarchy can be deduced from the tree structure of

$$\Gamma = \{(\hat{\xi}_s^{(x,t)}, s) : 0 \leq s \leq T\},$$

and we refer to [18], Section 2, and [11], Section 3, for a detailed description. To have a rigorous definition of the hierarchy of ancestors, we introduce a function $\phi_s$ that, at any time $s \geq 0$, maps the dual process into the set of sequences with integer values $S$ equipped with the usual lexicographic order $\ll$. We recall that for any $v, w \in S$,

$v \ll w$     if and only if     $v_k = w_k$     for $k = 1, 2, \ldots, n - 1$    and    $v_n < w_n$

for some integer $n \geq 1$. To facilitate the writing of $\phi_s$, we identify any $v \in S$ such that
$$v_n \neq 0 \quad \text{and} \quad v_k = 0 \quad \text{for } k \geq n+1$$
with the vector $(v_1, v_2, \ldots, v_n)$. The function $\phi_s$ is defined inductively in the following way. First, we go down the graphical representation starting at $(x, T)$, denote by $T - s_x$ the first time a death mark is encountered at site $x$ and let $\phi_s(x) = 0$ for any $0 \leq s \leq s_x$. Assume that a site $z$ is added to the dual process at time $T - u$ and denote by $T - s_z$ the first time we encounter a death mark at site $z$ by going backward in time starting from $(z, T - u)$. We then go back up from $(z, T - s_z)$ to $(z, T - u)$ and, each time we encounter the tip of an arrow, denote its space–time location by $(z_k, T - s_k)$, with $s_z > s_1 > \cdots > s_n > u$. For $k = 1, 2, \ldots, n$, denote by $T - s_{z_k}$ the first time a death mark is encountered at site $z_k$. We then set
$$\phi_s(z_k) = (u, k) \qquad \text{where } u = \phi_{s_k}(z) \qquad \text{for } s_k \leq s \leq s_{z_k}.$$
See Figure 3, left-hand picture. This naturally introduces a hierarchy in the set of ancestors where, given two sites $y$ and $z$ belonging to $\hat{\xi}_s^{(x,T)}$, site $y$ comes before site $z$ in the hierarchy if and only if $\phi_s(y) \ll \phi_s(z)$. The spatial location of the first ancestor at time $t - s$ is denoted by $\hat{\xi}_s^{(x,T)}(1)$. Later on, this ancestor will be called the *distinguished particle*.

The space–time point $(x, T)$ is said to *live forever* if $\hat{\xi}_s^{(x,T)} \neq \varnothing$ for all $s \geq 0$. Provided $(x, T)$ lives forever, the distinguished particle has a nice feature, namely that its path can be broken up at *renewal points* into independent and identically distributed pieces. Renewal points are defined as follows. If at time $s$ the distinguished particle jumps to a site $y$ such that $(y, s)$ lives forever, then $(y, s)$ is a renewal point. Let $(S_k, T_k)$ denote the $k$th renewal point, $X_i$ the spatial displacement between consecutive renewal points and $\tau_i$ the corresponding temporal displacement, that is,
$$S_k = x + \sum_{i=1}^{k} X_i \quad \text{and} \quad T_k = \sum_{i=1}^{k} \tau_i.$$
See Figure 3, right-hand picture.

PROPOSITION 2.1 ([18]). *Conditioned on the event that $(x, T)$ lives forever, the random vectors $\{(X_i, \tau_i)\}_{i \geq 1}$ are independent and identically distributed. Moreover, the tail distributions of $X_i$ and $\tau_i$ have exponential bounds, that is,*
$$P(\|X_i\| > t) \leq C e^{-\gamma t} \quad \text{and} \quad P(\tau_i > t) \leq C e^{-\gamma t}$$
*for appropriate $C < \infty$ and $\gamma > 0$.*

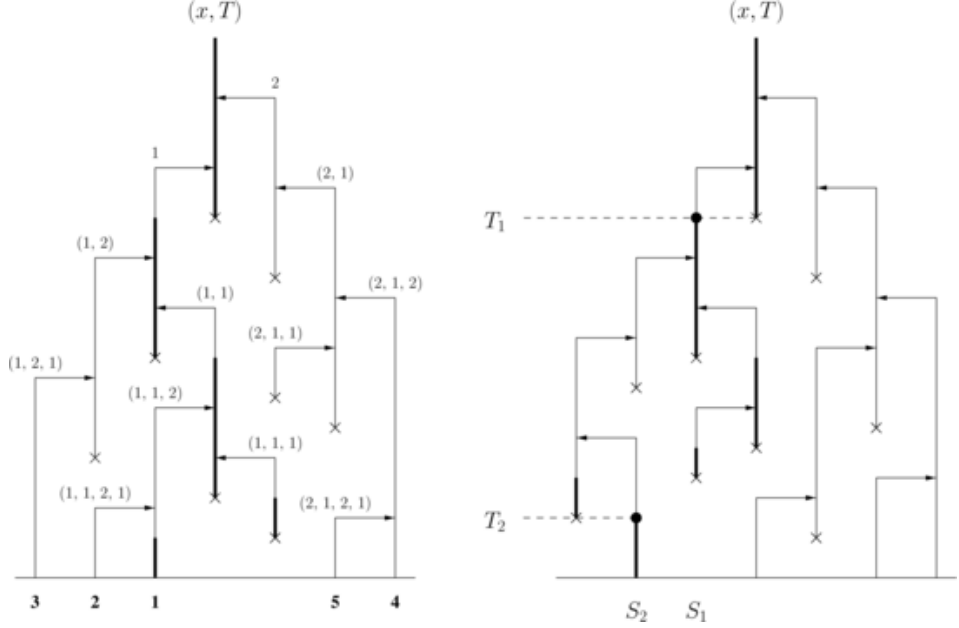

Fig. 3. *Tree structure of the dual process. Left-hand picture: Construction of the function $\phi_s$ and description of the ancestor hierarchy at time* 0. *Right-hand picture: Spatial and temporal locations of the renewal points. In both pictures, the bold lines refer to the path of the distinguished particle.*

Before proving Theorem 1, we need to describe an algorithm that allows us to deduce the type of $(x, T)$ from the (topological) structure of $\Gamma$ and the hierarchy of the ancestors. We look for the first ancestor in the hierarchy that does not land on an empty site at time 0. If this ancestor lands on a specialist (resp. generalist) at time 0, then $(x, T)$ will be occupied by a specialist (resp. generalist) of the same type unless the path crosses a **g**-arrow (resp. **s**-arrow) on its way up to $(x, T)$. If it lands on a specialist (resp. generalist) and the path crosses a **g**-arrow (resp. **s**-arrow), this specialist (resp. generalist) can block some other ancestors from determining the type of $(x, T)$ so we will need to remove these ancestors from the hierarchy. To do this, we follow the path on its way up to $(x, T)$ until we first encounter a **g**-arrow (resp. **s**-arrow), then remove from the hierarchy all the ancestors of the dual process starting at the point where this arrow is attached, repeat the same algorithm with the first remaining ancestor of the hierarchy at time 0, and so on. The algorithm halts whenever we find an ancestor that successfully reaches $(x, T)$. If none of the ancestors succeed, $(x, T)$ will be empty.

*Proof of Theorem* 1. It remains to prove Theorem 1 for $\alpha = \beta$. The proof is similar to the proof of Theorem 1 in [18]. We begin by assuming that the dual process starting at $(x,T)$ lives forever, an event with positive probability since $\alpha > \lambda_c$. We will prove that $P(\xi_T(x) = 3) \to 1$ as $T \to \infty$ on the set where $(x,T)$ lives forever. Without loss of generality, we can suppose that $x \in \mathcal{H}_1$. Proposition 2.1, together with the irreducibility of the random walk associated with the renewal points, implies that there is a subscript $k$ almost surely finite such that $S_k \in \mathcal{H}_2$. In particular, since what happens before and after each renewal is independent, we can inductively extract a subsequence of renewal points as follows. We set $k_0 = 0$ and, for any $i \geq 1$, let

$$k_i = \begin{cases} \min\{k > k_{i-1} : S_k \in \mathcal{H}_2\}, & \text{if } i \text{ is odd,} \\ \min\{k > k_{i-1} : S_k \in \mathcal{H}_1\}, & \text{if } i \text{ is even.} \end{cases}$$

By continuity of the path, the distinguished particle will cross, on its way up to $(x,T)$, at least one arrow that crosses from $\mathcal{H}_1$ to $\mathcal{H}_2$ (or vice versa) between time $T - T_{k_i}$ and time $T - T_{k_{i-1}}$. It follows from our graphical representation that such an arrow is a **g**-arrow which will block the specialists from determining the type of $(x,T)$.

To prove that $P(\xi_T(x) = 3) \to 1$ as $T \to \infty$ [provided $(x,T)$ lives forever], we will construct a sequence of ancestors that are good candidates to bring a generalist to $(x,T)$. We start the dual process at $(x,T)$ and determine the ordered set of ancestors after $T$ units of time by going backward in time. The first member of the sequence, denoted by $\eta_T(1)$, is the distinguished particle. We follow the path this ancestor takes to determine the type of $(x,T)$ by going forward in time until the first time we cross a **g**-arrow. We then discard all the ancestors that cross this **g**-arrow on their way up to $(x,T)$. The first ancestor that is left after discarding these ancestors, denoted by $\eta_T(2)$, is the second member of our sequence. Since, before discarding, the tree starting at each renewal point is linearly growing in time (in the sense described by the Shape Theorem in [17], page 128), we can find, by choosing $t$ sufficiently large, an integer $i \geq 1$ such that, even after discarding, the dual starting at $(S_{k_i}, T_{k_i})$ lives forever. This tells us that the path the second member will take to determine the type of $(x,T)$ goes through the $k_i$th renewal point. The definition of the sequence of times $T_{k_i}$, $i \geq 1$, implies that the second member crosses at least $k_i$ **g**-arrows on its way up to $(x,T)$. As previously, we discard all the ancestors that cross the first **g**-arrow encountered by the second member on its way up to $(x,T)$, define the third member of the sequence as the first remaining ancestor, and so on, until we run out of ancestors. Since different ancestors can occupy the same site, we now extract a further subsequence so that, at time 0, all the candidates are different. We start with $\eta_T(1)$ and discard all the members of $\{\eta_T(k) : k \geq 2\}$, that occupy the same site as $\eta_T(1)$. We then take the next ancestor that is left, and so on.

We now denote by $\eta_T$ the set of members of the subsequence we have just defined. With the possible exception which can occur if it is the last member of the subsequence [since it may happen that this member does not cross any **g**-arrow on its way up to $(x,T)$], if a site $z \in \eta_T$ is occupied by a specialist at time 0, then this specialist will be blocked on its way up to $(x,T)$ by a **g**-arrow. On the other hand, if $z$ is occupied by a generalist, this generalist will be able to bring its offspring up to $(x,T)$ since there are no **s**-arrows in our graphical representation (recall $\alpha = \beta$). In particular, if we denote by $\Theta_3$ the set of sites occupied at time 0 by a particle of type 3, we only need to prove that

$$P(\eta_T \cap \Theta_3 \neq \varnothing) \to 1 \quad \text{as } T \to \infty.$$

Since the tree growing out of $(x,T)$ is linearly growing in time, for any $\varepsilon > 0$ and $M > 0$, there exists a time $t_0 > 0$ such that

$$P(\operatorname{card} \eta_T < M) \leq \varepsilon \quad \forall T \geq t_0.$$

To conclude, we apply Lemma 9.14 from [14] which tells us that if $\xi_0$ is translation-invariant and $P(\xi_0(x) = 3) > 0$ then, given $\varepsilon > 0$, there is an $M > 0$ such that if $\operatorname{card} \eta_T \geq M$, then

$$P(\eta_T \cap \Theta_3 = \varnothing) \leq \varepsilon.$$

This, together with the previous inequality, completes the proof of Theorem 1.

**3. Proof of Theorem 2.** This section is devoted to the proof of Theorem 2, which describes the behavior of the stochastic process when $\alpha > \beta > \lambda_c$ and the space scale $L$ is large. To avoid cumbersome notation, we will prove the result only when $d = 2$, but would like to point out that when $d \geq 3$, each step of the proof is the same. The proof relies on a rescaling argument (see, e.g., [4, 8]). The basic idea is to show that, for any $\varepsilon > 0$, members of the family of processes under consideration, when viewed on suitable length and time scales, dominate a one-dependent, oriented site percolation process in which sites are open with probability at least $1 - \varepsilon$.

Let $B = \mathcal{H}_L = [-L, L]^2$. To prepare for a multiscale argument, we define a sequence of spatial boxes $J_n = [-L/n, L/n]^2$, $n = 1, 2, \ldots$. For any $z = (z_1, z_2)$ in $\mathbb{Z}^2$ and any positive integer $n \in \mathbb{N}^*$, we set

$$\Phi(z) = (2Lz_1, 2Lz_2), \qquad B(z) = \Phi(z) + B \quad \text{and} \quad J_n(z) = \Phi(z) + J_n.$$

Moreover, we tile $J_4(z)$ with $L^{0.1} \times L^{0.1}$ squares by setting

$$\pi(w) = (L^{0.1} w_1, L^{0.1} w_2), \qquad D = [-L^{0.1}/2, L^{0.1}/2]^2,$$
$$D(w) = \pi(w) + D, \qquad I_z = \{w \in \mathbb{Z}^2 : D(w) \subseteq J_4(z)\}.$$

To prove that specialists and generalists coexist, we start by introducing the following two collections of good events. We set $T = L^{1.5}$.

1. For $z_1$ and $z_2$ both even or both odd, we will say that $J_4(z)$ is **s**-*good* if $J_4(z)$ is void of 3's and has at least one particle of type 1 in each of the squares $D(w)$, $w \in I_z$. For $z_1$ and $z_2$ both even for even $k$ and $z_1$ and $z_2$ both odd for odd $k$, we will say that $(z,k)$ is **s**-*occupied* if the spatial box $J_4(z)$ is **s**-good at time $kT$.
2. For $z_1$ even and $z_2$ odd or $z_1$ odd and $z_2$ even, we will say that $J_4(z)$ is **g**-*good* if $J_4(z)$ is void of 1's and has at least one 3 in each of the squares $D(w)$, $w \in I_z$. For $z_1$ even and $z_2$ odd for even $k$ and $z_1$ odd and $z_2$ even for odd $k$, we will say that $(z,k)$ is **g**-*occupied* if the spatial box $J_4(z)$ is **g**-good at time $kT$.

The main ingredient needed to prove that the specialists survive is the following:

PROPOSITION 3.1. *Assume that $\alpha > \beta > \lambda_c$ and $J_4$ is **s**-good at time 0. For any $\varepsilon > 0$, the space scale parameter $L$ can then be chosen in such a way that*

$$P(J_4(1,1) \text{ is \textbf{s}-good at time } T) > 1 - \varepsilon.$$

The idea is that if some square $B(z) \subseteq \mathcal{H}_1$ is occupied by a cluster of specialists at a given time, then this cluster will invade the four diagonal squares with probability close to 1, and this will occur in less than $T$ units of time. It then follows from well-known results about oriented site percolation (see, e.g., [7, 10]) that the **s**-occupied sites percolate provided $\varepsilon > 0$ is sufficiently small, implying that the specialists survive. To prove that both consumers coexist, we still need to show that the **g**-occupied sites percolate as well. More precisely, we have, analogously:

PROPOSITION 3.2. *Assume that $\beta > \lambda_c$ and $J_4(0,1)$ is **g**-good at time 0. For any $\varepsilon > 0$, the space scale parameter $L$ can then be chosen in such a way that*

$$P(J_4(1,0) \text{ is \textbf{g}-good at time } T) > 1 - \varepsilon.$$

We illustrate the strategy of Proposition 3.1 in Figure 4. Let $z = (1,1)$ and $x \in J_4(z)$. We will then construct a dual path $A_t$ starting at $A_0 = (x, T)$ that lands at time 0 on the target set $J_4$ in the heart of the cluster of specialists. This selected path will have to block 3's, but not 1's, from determining the type of $(x, T)$. In particular, the important step will be to move $A_t$ from $J_4(z)$ to $J_4$ without escaping from the set $\mathcal{H}_1$, by crossing the corner at the point $Lz = (L, L)$. The proof of Proposition 3.2 is similar.

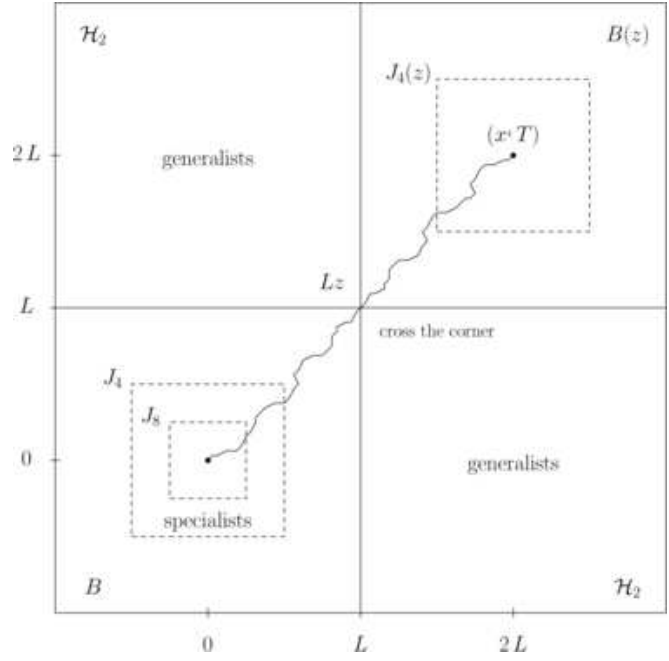

Fig. 4. *Picture of the selected path.*

*Construction of the selected path.* The construction of the selected path $A_t$ relies on the *repositioning algorithm* described in [11], supplemented with a restart argument. To define $A_t$, we will first construct a sequence of auxiliary processes for $k = 0, 1, \ldots$, denoted by $\{\hat{\xi}_t^k\}_{t \geq 0}$, called the $k$th *test path*, with $\hat{\xi}_t^k$ taking values in $\mathbb{Z}^2$, and then show in a series of lemmas that we can use the test paths to construct $A_t$. We will say that a renewal point is *associated with an* **s**-*arrow* if the first arrow a particle crosses, starting at this renewal point and moving up the graphical representation, is an **s**-arrow (see Figure 5, left-hand picture). Furthermore, let $\Delta$ denote the straight line going through the origin $(0,0)$ and the point $Lz$, and denote by $\text{dist}(\cdot, \cdot)$ the Euclidean distance. We also let $\delta > 0$ be a constant, to be fixed later.

The test paths are defined inductively (see Figure 5, right-hand picture). To construct the first test path $\{\hat{\xi}_t^1\}_{t \geq 0}$, we define a sequence of stopping times $\sigma_{1,j}$, $j = 0, 1, 2, \ldots$, with $\sigma_{1,0} = 0$. The process $\hat{\xi}_t^1$ starts at $\hat{\xi}_0^1 = x$ and follows the path of the first ancestor starting at $(x, T)$ until the first time, denoted by $\sigma_{1,1}$, that it jumps to a renewal point associated with an **s**-arrow. At time $\sigma_{1,1}$, we either leave $\hat{\xi}_t^1$ where it is at that time or we reposition $\hat{\xi}_t^1$. To determine whether and where to reposition $\hat{\xi}_t^1$, we denote the location of the second ancestor in the hierarchy at time $\sigma_{1,1}$, if it exists, by $B_{\sigma_{1,1}}$. If

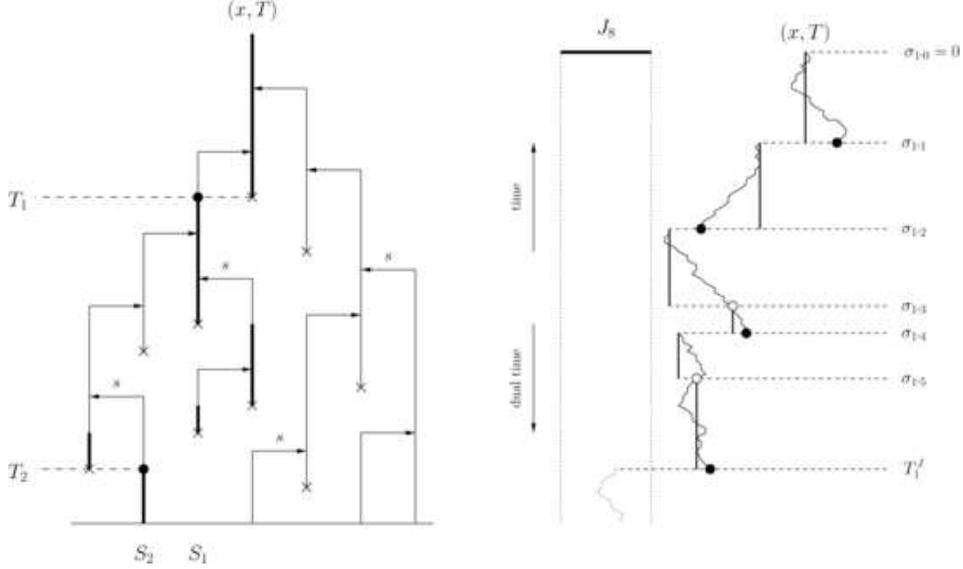

FIG. 5. *Left-hand picture: Picture of the renewal points. The first renewal point $(S_1, T_1)$ is not associated with an **s**-arrow. The second one, namely $(S_2, T_2)$, is associated with an **s**-arrow. In particular, $\sigma_{1,1} = T_2$. Right-hand picture: Picture of the first test path (solid path). The bold straight lines refer to the jump process $S_t^1$. The points represent the renewal points associated with an **s**-arrow. The test path $\hat{\xi}_t^1$ is repositioned at closed points and left where it was at that time at open points.*

$B_{\sigma_{1,1}}$ exists and the dual process starting at $(B_{\sigma_{1,1}}, T - \sigma_{1,1})$ lives forever, then we have the following three alternatives:

1. If $\mathrm{dist}(\hat{\xi}^1_{\sigma_{1,1}}, \Delta) > \delta$ and $\mathrm{dist}(B_{\sigma_{1,1}}, \Delta) < \mathrm{dist}(\hat{\xi}^1_{\sigma_{1,1}}, \Delta)$, then we set $\hat{\xi}^1_{\sigma_{1,1}+} = B_{\sigma_{1,1}}$.
2. If $\mathrm{dist}(\hat{\xi}^1_{\sigma_{1,1}}, \Delta) \leq \delta$ and $\mathrm{dist}(B_{\sigma_{1,1}}, 0) < \mathrm{dist}(\hat{\xi}^1_{\sigma_{1,1}}, 0)$, then we set $\hat{\xi}^1_{\sigma_{1,1}+} = B_{\sigma_{1,1}}$.
3. Otherwise, we set $\hat{\xi}^1_{\sigma_{1,1}+} = \hat{\xi}^1_{\sigma_{1,1}}$.

If $B_{\sigma_{1,1}}$ does not exist, or the dual process starting at $(B_{\sigma_{1,1}}, T - \sigma_{1,1})$ does not live forever, then we set $\hat{\xi}^1_{\sigma_{1,1}+} = \hat{\xi}^1_{\sigma_{1,1}}$.

To define the next stopping time $\sigma_{1,2}$, we start a new dual process at $(\hat{\xi}^1_{\sigma_{1,1}+}, T - \sigma_{1,1})$ and follow the path of the first ancestor until the first time it jumps to a renewal point that is associated with an **s**-arrow. Denote this time by $\sigma_{1,2}$. [Note that we followed the path of the first ancestor of the dual process starting at $(\hat{\xi}^1_{\sigma_{1,1}+}, T - \sigma_{1,1})$ for $\sigma_{1,2} - \sigma_{1,1}$ units of time, as shown in Figure 5, right-hand picture.] We apply the repositioning algorithm again as described above.

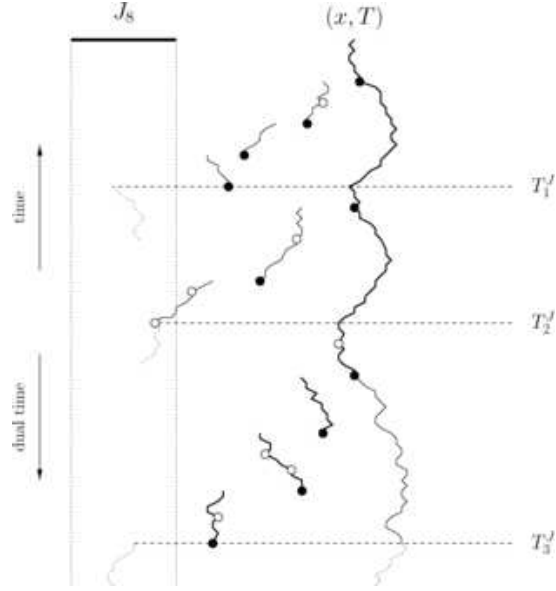

FIG. 6. *Picture of the first three test paths. The selected path $A_t$, in bold lines, is equal to the third test path, which is the first successful path. As in Figure 5, right-hand picture, the points represent the renewal points associated with an **s**-arrow. The repositioning algorithm is applied successfully at filled points.*

Into the first test path, we embed a jump process $S_t^1$, $t \geq 0$, as follows. The process $S_t^1$ stays put except at the stopping times $\sigma_{1,j}$, $j = 0, 1, \ldots$, that is, $S_t^1 = \hat{\xi}_{\sigma_{1,j}+}^1$ for $\sigma_{1,j} \leq t < \sigma_{1,j+1}$. The effect of the repositioning algorithm is that, whenever possible, we move the jump process $S_t^1$ closer to 0 while staying within a finite distance $\delta$ of $\Delta$.

We set $T_0^J = 0$ and let $T_1^J$ denote the first time $S_t^1$ jumps inside $J_8$. We call $\hat{\xi}_t^1$ a *successful path* if a particle, starting at $(\hat{\xi}_{T_1^J}^1, T - T_1^J)$ and going up to $(x, T)$ by taking the path we have just defined, does not cross any **g**-arrow. If the 1st path is not successful, we define a 2nd path $\hat{\xi}_t^2$ as follows. For the first $T_1^J$ units of time, the test path follows the path of the distinguished particle starting at $\hat{\xi}_0^2 = x$, and from time $T_1^J$ on moves according to the repositioning algorithm. Into this process, we embed the jump process $S_t^2$ and denote by $T_2^J$ the first time $S_t^2$ jumps inside $J_8$. We repeat the procedure until the first successful path is found, and call this path $A_t$ (see Figure 6).

*Existence of a successful path.* Since the $k$th test path $\hat{\xi}_t^k$ follows the distinguished particle starting at $(x, T)$ from time 0 to time $T_{k-1}^J$, the first step is to prove that, with probability close to 1, the first ancestor $\hat{\xi}_t^{(x,T)}(1)$ does not escape from the box $J_2(z)$ by time $T = L^{1.5}$. This allows us to locate the path when we start applying the repositioning algorithm. It also implies that our $k$th path does not cross any **g**-arrows between time 0 and time $T_{k-1}^J$ and is therefore a successful path.

LEMMA 3.3. *Assume that $\alpha > \beta > \lambda_c$. If $x \in J_4(z)$ and $T = L^{1.5}$, then*
$$P(\hat{\xi}_t^{(x,T)}(1) \notin J_2(z) \text{ for } t \in [0,T]) \leq C_1 \exp(-\gamma_1 L^{0.2})$$
*for appropriate $C_1 < \infty$ and $\gamma_1 > 0$.*

PROOF. Since our graphical representation is invariant under translation of the vector $2Lz$, we can assume that $x \in J_4$ and work with $J_2$ instead of $J_2(z)$. We will prove that, with probability close to 1, (i) the renewal points do not escape from the box $x + J_8$ and (ii) between two renewal points, the spatial displacements of the distinguished particle do not exceed $L/8$, an event that we define as
$$G_k = \{|S_k - \hat{\xi}_t^{(x,T)}(1)| \leq L/8 \text{ for } T_k \leq t < T_{k+1}\}.$$
Let $m_1 = \mathbb{E}\tau_1$ and denote by $S_t$ the location of the renewal point $S_n$ at time $t$, and by $N_T$ the number of renewal points between time 0 and time $T$. We then have
$$P(\hat{\xi}_t^{(x,T)}(1) \notin J_2 \text{ for some } t \leq T)$$
$$\leq P(N_T > 2T/m_1)$$
$$\quad + P(G_k^c \text{ for some } 0 \leq k \leq N_T; N_T \leq 2T/m_1)$$
$$\quad + P(S_t - x \notin J_8 \text{ for some } 0 \leq t \leq T; N_T \leq 2T/m_1).$$
By using large deviation estimates, we get
$$P(N_T > 2T/m_1) \leq C_2 \exp(-\gamma_2 T)$$
for appropriate $C_2 < \infty$ and $\gamma_2 > 0$. Furthermore, the proof of Proposition 2.1 implies that
$$P(G_k^c \text{ for some } 0 \leq k \leq N_T; N_T \leq 2T/m_1)$$
$$\leq P(G_k^c \text{ for some } 1 \leq k \leq 2T/m_1)$$
$$\leq 2T/m_1 P(G_1^c) \leq 2T/m_1 C \exp(-\gamma L/8).$$
Finally, by letting $S_t^{(1)}$ denote the 1st coordinate of $S_t$ and using the reflection principle, we get
$$P(S_t - x \notin J_8 \text{ for some } 0 \leq t \leq T; N_T \leq 2T/m_1)$$
$$\leq 2P(|S_t^{(1)} - S_0^{(1)}| \geq L/8 \text{ for some } 0 \leq t \leq T; N_T \leq 2T/m_1)$$
$$\leq 4P(|S_{2T/m_1}^{(1)} - S_0^{(1)}| \geq L/8).$$

Now, Chebyshev's inequality implies that, for any $\theta > 0$,

$$P(|S^{(1)}_{2T/m_1} - S^{(1)}_0| \geq L/8) \leq e^{-\theta L/8} \prod_{i=1}^{2T/m_1} \mathbb{E} \exp(\theta X_i^{(1)})$$

$$\leq \exp\left(-\frac{\theta L}{8} + \frac{2L^{1.5}}{m_1} \log \phi(\theta)\right),$$

where $\phi(\theta)$ is the moment-generating function of $X_1^{(1)}$. Since $\mathbb{E} X_1^{(1)} = 0$ and $\operatorname{Var} X_1^{(1)} < \infty$, we have $\log \phi(\theta) \leq C_3 \theta^2$ for some $C_3 < \infty$ and for sufficiently small $\theta$. By taking $\theta = L^{-0.8}$ in the last expression, we can conclude that

$$P(S_t - x \notin J_8 \text{ for some } 0 \leq t \leq T; N_T \leq 2T/m_1) \leq 4 \exp(-L^{0.2}/10)$$

for sufficiently large $L$. Putting things together, we finally get

$$P(\hat{\xi}_t^{(x,T)}(1) \notin J_2 \text{ for some } t \leq T)$$
$$\leq C_2 \exp(-\gamma_2 L^{1.5}) + 2L^{1.5}/m_1 C \exp(-\gamma L/8) + 4\exp(-L^{0.2}/10).$$

This completes the proof of the lemma. □

The second step in proving the existence of a successful path is to investigate the behavior of the process $\hat{\xi}_t^m$ between time $T^J_{m-1}$ and time $T^J_m$, that is, when the repositioning algorithm is applied. The aim is to prove that there is a positive probability that, between these times, the $m$th test path does not cross any **g**-arrows on its way to $J_8$. Since this will only happen if the path crosses the corner as shown in Figure 4, we will call this event *to successfully cross the corner*. Let $\ell$ be a positive integer to be fixed later and, for any $k \geq 0$, set

$$a_k = Lz + 2^k(\ell, 0), \qquad b_k = Lz + 2^k(0, \ell),$$
$$c_k = Lz + 2^k(\ell, 2\ell), \qquad d_k = Lz + 2^k(2\ell, \ell).$$

Moreover, we consider the sequence of segments

$$I_k = [3a_k/8 + 5b_k/8, 5a_k/8 + 3b_k/8],$$

as well as the rapidly increasing sequence of squares given by $D_k = (a_k, b_k, c_k, d_k)$. That is, $I_k$ is the segment with endpoints $3a_k/8 + 5b_k/8$ and $5a_k/8 + 3b_k/8$, and $D_k$ the square with corners $a_k, b_k, c_k$ and $d_k$. We can assume that $\ell$ and $L$ are chosen so that there exists $K^\star$ with $L = 2^{K^\star}\ell$. $I_{K^\star}$ is then contained in the diagonal line connecting the points $Lz + (2L, 0)$ and $Lz + (0, 2L)$ (see Figure 7 for a picture). Assume that the $m$th test path starts at $x \in I_{k+1}$. We will prove that, with a positive probability (uniform in $k$), we can bring our test path from $x$ to $I_0$ without escaping from the decreasing sequence of squares $D_k, D_{k-1}, \ldots, D_0$. We say that the repositioning algorithm was *applied successfully* if the new path is chosen. Let $M_t$ denote the number of times the repositioning algorithm has been applied successfully by time $t$.

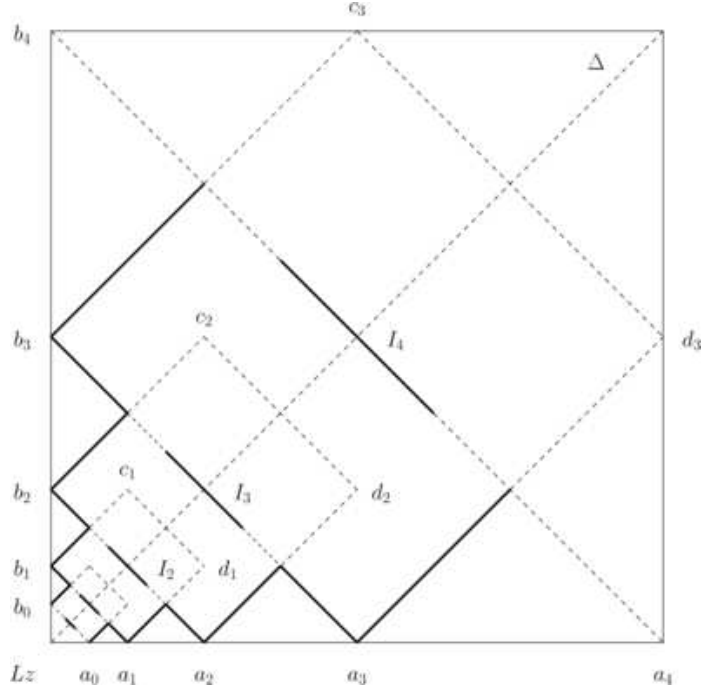

Fig. 7. *The sequence of squares.*

LEMMA 3.4. *Assume that $\alpha > \beta > \lambda_c$. There exists $\gamma_4 > 0$ so that*
$$P(M_t \leq \gamma_4 t) \leq C_5 \exp(-\gamma_5 t)$$
*for appropriate $C_5 < \infty$ and $\gamma_5 > 0$.*

PROOF. This is the same as the proof of Lemma 3.4 in [11]. □

In Lemma 3.5, we will show that, with high probability, a test path can move from $I_1$ to $I_0$ without leaving $D_0$. A corollary to this lemma is that, with probability close to 1, a test path can move from $I_k$ to $I_{k-1}$ without leaving $D_{k-1}$. Applying this lemma repeatedly will then show that, with positive probability, a test path can successfully cross the corner starting anywhere in $J_2(z)$.

LEMMA 3.5. *If $\alpha > \beta > \lambda_c$ and the mth test path starts within distance $\ell/8$ of $I_1$, then*

$P($the test path escapes from $D_0$ before getting within $\ell/16$ of $I_0)$
$$\leq C_6 \exp(-\gamma_6 \ell^{0.4})$$
*for sufficiently large $\ell$ and appropriate $C_6 < \infty$ and $\gamma_6 > 0$.*

PROOF. To lighten the notation, we will write $\sigma_k$ instead of $\sigma_{m,k}$. Let $\theta_1$ (resp. $\theta_2$) denote the first time $\text{dist}(S_t^m, \Delta) \leq \delta$ [resp. $\hat{\xi}_t^m$ crosses either segment $(a_0 d_0)$ or $(b_0 c_0)$], and $\pi_\Delta$ denote the orthogonal projection on the straight line $\Delta$. It follows from the repositioning algorithm that there exists $m_2 > 0$ such that, on the set where $\text{dist}(S_{\sigma_{k-1}}^m, \Delta) > \delta$,

$$\mathbb{E}[\text{dist}(S_{\sigma_k}^m, \Delta) - \text{dist}(S_{\sigma_{k-1}}^m, \Delta) | \mathcal{F}_{\sigma_{k-1}}] \leq -m_2$$

where $\mathcal{F}_{\sigma_{k-1}} = \sigma(S_{\sigma_0}^m, \ldots, S_{\sigma_{k-1}}^m)$. This, together with Lemma 3.4, implies that we can find constants $C_7$, $C_8$, $\gamma_8 \in (0, \infty)$ such that

$$P(\theta_1 > C_7 \ell \text{ or } \theta_1 > \theta_2) \leq C_8 \exp(-\gamma_8 \ell).$$

Furthermore, using the notation introduced in Lemma 3.3, we get

$$P(\text{dist}(\hat{\xi}_t^m, (a_1 b_1)) > \ell/2 \text{ for some } t \leq \theta_1)$$
$$\leq P(\theta_1 > C_7 \ell \text{ or } \theta_1 > \theta_2) + P(N_{\theta_1} > 2C_7 \ell/m_1; \theta_1 \leq C_7 \ell)$$
$$+ P(\text{the largest spatial displacement is greater than } \ell/16;$$
$$N_{\theta_1} \leq 2C_7 \ell/m_1)$$
$$+ P(\text{dist}(S_t^m, (a_1 b_1)) > \ell/4 \text{ for some } t \leq \theta_1; N_{\theta_1} \leq 2C_7 \ell/m_1).$$

Since $\pi_\Delta(S_{\sigma_k}^m) - \pi_\Delta(S_{\sigma_{k-1}}^m)$ has mean 0 (and finite variance) before time $\theta_1$, the arguments of the proof of Lemma 3.3 (with $\theta = \ell^{-0.6}$) together with the exponential bound on $\theta_1$ imply that

$$P(\text{dist}(\hat{\xi}_t^m, (a_1 b_1)) > \ell/2 \text{ for some } t \leq \theta_1)$$
$$\leq C_8 \exp(-\gamma_8 \ell) + C_2 \exp(-\gamma_2 C_7 \ell)$$
$$+ 2C_7 \ell/m_1 C \exp(-\gamma \ell/16) + 2 \exp(-\ell^{0.4}/8).$$

Now that $\hat{\xi}_t^m$ is close to $\Delta$, the next step is to bring it to the target region $I_0$ without escaping from the square $D_0$. Let $\theta_3 > \theta_1$ (resp. $\theta_4 > \theta_1$) denote the first time $S_t^m$ crosses the straight line $(a_0 b_0)$ [resp. $(c_0 d_0)$], and let $K$ denote the number of times the repositioning algorithm has been applied between time $\theta_1$ and time $\theta_3$. First of all, since $S_t^m$ has a drift toward $\Delta$ when $\text{dist}(S_{\sigma_{k-1}}^m, \Delta) > \delta$, there is a positive constant $\rho$ such that

$$P(\text{card}(\sigma_k \in [\theta_1, \theta_3] : \text{dist}(S_{\sigma_k}^m, \Delta) \leq \delta) < \rho K) \leq C_9 \exp(-\gamma_9 \ell)$$

for appropriate $C_9 < \infty$ and $\gamma_9 > 0$. That is, with probability close to 1, the fraction of renewal points associated with **s**-arrows that are within distance $\delta$ from $\Delta$ is at least $\rho$. Furthermore, there exists $m_3 > 0$ such that, on the set where $\text{dist}(S_{\sigma_{k-1}}^m, \Delta) \leq \delta$,

$$\mathbb{E}[\text{dist}(S_{\sigma_k}^m, 0) - \text{dist}(S_{\sigma_{k-1}}^m, 0) | \mathcal{F}_{\sigma_{k-1}}] \leq -m_3$$

which, together with the previous estimate, implies that there is a $C_{10} < \infty$ such that

$$P(\theta_3 - \theta_1 > C_{10}\ell \text{ or } \theta_3 > \theta_4) \leq C_{11} \exp(-\gamma_{11}\ell)$$

for appropriate $C_{11} < \infty$ and $\gamma_{11} > 0$. The same decomposition as above then implies that

$$P(\text{dist}(\hat{\xi}_t^m, \Delta) > \ell/8 \text{ for some } \theta_1 \leq t \leq \theta_3) \leq C_{12} \exp(-\gamma_{12}\ell)$$

for appropriate $C_{12} < \infty$ and $\gamma_{12} > 0$. In conclusion, there exist $C_6, \gamma_6 \in (0, \infty)$ such that

$P(\text{the test path escapes from } D_0 \text{ before getting within } \ell/16 \text{ of } I_0)$
$\leq P(\text{dist}(\hat{\xi}_t^m, (a_1 b_1)) > \ell/2 \text{ for some } t \leq \theta_1) + P(\theta_1 > C_7\ell \text{ or } \theta_1 > \theta_2)$
$\quad + P(\text{dist}(\hat{\xi}_t^m, \Delta) > \ell/8 \text{ for some } \theta_1 \leq t \leq \theta_3)$
$\quad + P(\theta_3 - \theta_1 > C_{10}\ell \text{ or } \theta_3 > \theta_4)$
$\leq C_6 \exp(-\gamma_6 \ell^{0.4}).$

This completes the proof of the lemma. □

Relying on Lemma 3.5, we can now prove that there is a strictly positive probability that, starting at $x \in J_2(z)$, the $m$th test path reaches $J_8$ without crossing any **g**-arrow. This, together with independence of temporally separated paths, will imply that a geometric number of trials suffices to find a successful path.

LEMMA 3.6. *Assume that $\alpha > \beta > \lambda_c$. There exists $p_1 > 0$ such that*

$$P(\text{the } m\text{th test path successfully crosses the corner}) > p_1$$

*for all sufficiently large $L$.*

PROOF. This is a three-step process in which we will prove that, with positive probability, we can bring our path from $J_2(z)$ to $I_0$, then from $I_0$ to $I_0 - (\ell, \ell)$, and finally from $I_0 - (\ell, \ell)$ to $J_8$ without crossing any **g**-arrows (see Figures 4 and 7). To estimate our first good event, we apply Lemma 3.5 repeatedly to obtain

$P(\text{starting within } L/16 \text{ of } I_{K^\star},$
$\quad \text{the path crosses a } \mathbf{g}\text{-arrow before getting within } \ell/16 \text{ of } I_0)$
$\leq C_6 \sum_{k=0}^{\infty} \exp(-\gamma_6 (2^k \ell)^{0.4}) \leq 2 C_6 \exp(-\gamma_6 \ell^{0.4})$

for sufficiently large $\ell$. Note that the previous estimate holds regardless of the value of $L$. Moreover, it tells us that we can fix $\ell$ so that our first good event occurs with probability $p_2 > 0$. The same reasoning that led to the estimate of the first good event can be used to estimate the third good event, since it is symmetric to our first good event. We reverse time in Lemma 3.5 and find that if a test path starts within $\ell/16$ of $I_0$, then the probability that the test path escapes from $D_0$ before getting within $\ell/8$ of $I_1$ is less than $C_6 \exp(-\gamma_6 \ell^{0.4})$. Iterating as before, it follows that, with high probability, we can move the particle inside $J_8$ without crossing any **g**-arrows. Now that $\ell$ is fixed, there is a positive probability $p_3(\ell) > 0$ that the $m$th test path goes from $I_0$ to $I_0 - (\ell, \ell)$ without crossing any **g**-arrow, since the interaction neighborhood we chose has (at least) the eight nearest neighbors. Finally, to be able to start from anywhere in $J_2(z)$, we use the repositioning algorithm to move the path first from $x \in J_2(z)$ toward $\Delta$ and then toward $I_{K^\star}$. Using similar estimates as in Lemma 3.5, it follows that

$P($the test path escapes from $J_1(z)$ before getting within $L/16$ of $I_{K^\star})$
$$\leq C_6 \exp(-\gamma_6 L^{0.4})$$

for sufficiently large $L$ and appropriate constants $C_6, \gamma_6 \in (0, \infty)$. In conclusion, the lemma follows by taking $p_1 = p_2^2 p_3(\ell) > 0$. □

With Lemmas 3.3 and 3.6 in hand, we are now ready to prove that, with high probability, the first successful path reaches $J_8$ by time $T - \sqrt{L}$. More precisely, we have the following:

LEMMA 3.7. *Let $T = L^{1.5}$ and $N = \min\{k \geq 1 : \hat{\xi}_t^k \text{ is a successful path}\}$. It then follows that*
$$P(T_N^J > T - \sqrt{L}) \leq C_{13} \exp(-\gamma_{13} L^{0.2})$$
*for $\alpha > \beta > \lambda_c$ and appropriate $C_{13} < \infty$ and $\gamma_{13} > 0$.*

PROOF. The first step is to prove that the temporal displacement $T_{k+1}^J - T_k^J$ is small enough such that, by taking sufficiently large $L$, the number of test paths run by time $T$ can be made arbitrarily large. The arguments used in the proof of Lemma 3.5 to estimate the stopping times $\theta_1$ and $\theta_3$ imply that there exists $C_{14} < \infty$ such that
$$P(T_{k+1}^J - T_k^J > C_{14} L) \leq C_{15} \exp(-\gamma_{15} L)$$
for appropriate $C_{15} < \infty$ and $\gamma_{15} > 0$. In particular, since the increments $T_{k+1}^J - T_k^J$ are independent, it remains to be shown that a geometrical number of trials suffices to find a successful path. Let $\Omega$ denote the event that the distinguished particle starting at $(x,T)$ does not cross any **g**-arrow by time $T$. From Lemma 3.3 we know that

$$P(\Omega) \geq 1 - C_1 \exp(-\gamma_1 L^{0.2})$$

and from Lemma 3.6 that

$P($the $k$th test path does not

cross any **g**-arrow by time $T_k^J | \Omega; T_{k-1}^J < T) > p_1,$

regardless of the values of $k$ and $L$. This implies that

$$\begin{aligned}
P(T_N^J > T - \sqrt{L}) &\leq P(T_N^J > T - \sqrt{L}; N > L^{0.4}) \\
&\quad + P(T_N > T - \sqrt{L}; N \leq L^{0.4}; \Omega) \\
&\quad + P(T_N > T - \sqrt{L}; N \leq L^{0.4}; \Omega^c) \\
&\leq (1-p_1)^{L^{0.4}} + L^{0.4} C_{15} \exp(-\gamma_{15} L) + C_1 \exp(-\gamma_1 L^{0.2}) \\
&\leq C_{13} \exp(-\gamma_{13} L^{0.2})
\end{aligned}$$

for appropriate $C_{13} < \infty$ and $\gamma_{13} > 0$. The lemma follows. □

*Invading the spatial box* $J_4(z)$. By Lemma 3.7, our selected path $A_t$ reaches the target set $J_8$ by time $T - \sqrt{L}$ with probability close to 1. The last step, similar to the one described in [11], is to show that such a path gives us a good opportunity to bring a 1 up to $(x,T)$.

To define $A_t$ between time $T_N^J$ and time $\tau = T - \sqrt{L}$, we pretend that our selected path follows the path determined by the algorithm of the first ancestor starting at $(A_{T_N^J}, T - T_N^J)$. The next step is to prove that, between times $T_N^J$ and $\tau$, the selected path does not escape from $J_6$. We will then let the dual process spread out for the remaining $\sqrt{L}$ time units and sprinkle the target region $J_4$ where the 1's live.

LEMMA 3.8. *Assume that* $\alpha > \beta > \lambda_c$. *Let* $T = L^{1.5}$ *and* $\tau = T - \sqrt{L}$. *It then follows that*

$$P(A_t \notin J_6 \text{ for some } T_N^J \leq t \leq \tau) \leq C_{16} \exp(-\gamma_{16} L^{0.2})$$

*for appropriate* $C_{16} < \infty$ *and* $\gamma_{16} > 0$.

PROOF. This is, with some minor modifications, the same as the proof of Lemma 3.3. □

To complete the proof of Proposition 3.1, we still need to prove that $(A_\tau, \sqrt{L})$ will be occupied by a 1, provided $J_4$ is **s**-good at time 0. Since the selected path has been constructed to block 3's, this 1 will determine the type of $(x, T)$ unless another 1 succeeds earlier. This will prove that the specialists present in $J_4$ invade the neighboring box $J_4(z)$ with probability close to 1.

LEMMA 3.9. *Assume that $\alpha > \beta > \lambda_c$ and that $J_4$ is **s**-good at time 0. If $A_\tau \in J_6$, there are $C_{17} < \infty$ and $\gamma_{17} > 0$ such that*

$$P((A_\tau, \sqrt{L}) \text{ is not occupied by a } 1) \leq C_{17} \exp(-\gamma_{17} L^{0.1}).$$

PROOF. We run the dual process starting at $(A_\tau, \sqrt{L})$ for $\sqrt{L} - L^{0.2}$ units of time. It follows from the Shape Theorem (see [17], page 128) and the ergodicity of the upper invariant measure of the contact process (see Proposition 2.16, page 143, in [16]) that, except for a probability smaller than $C_{18} \exp(-\gamma_{18} L^{0.1})$, we can select $L^{0.1}$ sites at time $L^{0.2}$ which are contained in the dual process starting at $(A_\tau, \sqrt{L})$ such that (i) all these sites are contained in $J_5$, (ii) they are at least $L^{0.3}$ units apart from each other, and (iii) none of the duals starting at these sites uses the same part of the graphical representation for the remaining $L^{0.2}$ units of time. Since each of the duals has a positive probability of surviving, and the square $J_4$ is **s**-good at time 0, each of these $L^{0.1}$ sites has probability $\eta > 0$ of being occupied by a 1. In conclusion,

$$P((A_\tau, \sqrt{L}) \text{ is not occupied by a } 1) \leq C_{18} \exp(-\gamma_{18} L^{0.1}) + (1-\eta)^{L^{0.1}}.$$

This completes the proof of the lemma. □

The proof of Proposition 3.1 is now straightforward. First of all, by combining Lemmas 3.7, 3.8 and 3.9, it follows that if $J_4$ is **s**-good at time 0, then, for any $x \in J_4(z)$,

$$P(\xi_T(x) = 3) \leq C_{13} \exp(-\gamma_{13} L^{0.2}) + C_{16} \exp(-\gamma_{16} L^{0.2}) + C_{17} \exp(-\gamma_{17} L^{0.1})$$
$$\leq C_{19} L^{-3}$$

for some $C_{19} < \infty$. In particular, since there are $(L/2 + 1)^2$ sites in $J_4(z)$, it follows that

$$P(\xi_T(x) = 3 \text{ for some } x \in J_4(z)) \leq C_{19} L^{-1} \leq \varepsilon/2$$

for sufficiently large $L$. In other respects, the process dominates a one-color contact process with parameter $\beta > \lambda_c$, so the probability that there exists an $L^{0.1} \times L^{0.1}$ square $D(w) \subseteq J_4(z)$ that is empty can be bounded by $C_{20} L^{1.8} \exp(-\gamma_{20} L^{0.2}) \leq \varepsilon/2$. Proposition 3.1 then follows.

PROOF OF PROPOSITION 3.2. The proof of Proposition 3.2 is almost the same. The construction of the selected path, however, is slightly different, since there are no **s**-arrows in the squares $B(0,1)$ and $B(1,0)$. We will use the same notation to describe the construction of the new selected path.

Let $e_1 = (1,0)$ and $x \in J_4(e_1)$. As previously, $A_t$ will be the first successful path in a sequence of test paths. The first difference is that, since $B(e_1) \subseteq \mathcal{H}_2$, and the 1's cannot give birth through an arrow that points into $\mathcal{H}_2$, we can now apply the repositioning algorithm at each renewal point, regardless of whether it is associated with an **s**-arrow or not. On the other hand, the **s**-arrows are now prohibited, so we remove them from the graphical representation. We run the dual process, starting at $(x, T)$ for $T$ units of time, and break up the path of the distinguished particle at renewal points. Since $\beta > \lambda_c$, the event that the dual process lives forever still has a positive probability, even after removal of the **s**-arrows.

To define the 1st test path $\hat{\xi}_t^1$, we first set $\sigma_{1,0} = 0$. The process $\hat{\xi}_t^1$ starts at $\hat{\xi}_0^1 = x$ and follows the path of the first ancestor of $(x, T)$ until the first time, denoted by $\sigma_{1,1}$, that it jumps to a renewal point. As previously, the position of the test path at time $\sigma_{1,1}$ is determined by a repositioning algorithm. The process then follows the path described by the algorithm of the first ancestor until the next time, denoted by $\sigma_{1,2}$, that it jumps to a new renewal point, and so on. Let $B_{\sigma_{1,k}}$ denote the second ancestor, if it exists, of the previous renewal at time $\sigma_{1,k}$. Since $\beta > \lambda_c$, there is a positive probability that $(B_{\sigma_{1,k}}, T - \sigma_{1,k})$ lives forever. To describe the new repositioning algorithm, we define as above an embedded jump process $S_t^1$ which stays put except at times $\sigma_{1,k}$. Let $e_2 = (0,1)$ and let $\Delta'$ denote the straight line going through $Lz$ and $2Le_2$. The jump process $S_t^1$ evolves according to the following rules. If site $B_{\sigma_{1,k}}$ does not exist, or if the dual starting at $(B_{\sigma_{1,k}}, T - \sigma_{1,k})$ does not live forever, then we set $S_{\sigma_{1,k}}^1 = \hat{\xi}_{\sigma_{1,k}}^1$. Otherwise, we have the following three alternatives:

1. If $\text{dist}(\hat{\xi}_{\sigma_{1,k}}^1, \Delta') > \delta$ and $\text{dist}(B_{\sigma_{1,k}}, \Delta') < \text{dist}(\hat{\xi}_{\sigma_{1,k}}^1, \Delta')$, then we set $S_{\sigma_{1,k}}^1 = B_{\sigma_{1,k}}$.
2. If $\text{dist}(\hat{\xi}_{\sigma_{1,k}}^1, \Delta') \leq \delta$ and $\text{dist}(B_{\sigma_{1,k}}, 2Le_2) < \text{dist}(\hat{\xi}_{\sigma_{1,k}}^1, 2Le_2)$, then we set $S_{\sigma_{1,k}}^1 = B_{\sigma_{1,k}}$.
3. Otherwise, we set $S_{\sigma_{1,k}}^1 = \hat{\xi}_{\sigma_{1,k}}^1$.

In other words, we now move the jump process toward the target set $J_8(e_2)$ in a direction perpendicular to the straight line $\Delta$.

As previously, the 2nd test path follows the path of the first ancestor of $(x, T)$ until the 1st test path reaches $J_8(e_2)$, at which time it moves according to the repositioning algorithm, and so on. We now call $\hat{\xi}_t^k$ a successful path if a particle in $J_8(e_2)$ following this path by going forward in time does not cross any arrow that points to $\mathcal{H}_1$ on its way up to $(x, T)$. In particular, if this particle is of type 3, it will determine the type of $(x, T)$ unless another 3 succeeds earlier.

It is not difficult to see that, with this new construction, all our estimates still hold. In particular, if the square $J_4(e_2)$ is **g**-good at time 0, then

$$P(\xi_T(x) = 1 \text{ for some } x \in J_4(z)) \leq C_{21} L^{-1}$$

for a suitable $C_{21} < \infty$. Proposition 3.2 and Theorem 2 then follow.  $\square$

## REFERENCES

[1] ANDERSON, K. and NEUHAUSER, C. (2002). Patterns in spatial simulations—are they real? *Ecological Modelling* **155** 19–30.
[2] BERNAYS, E. A. (1989). The value of being a resource specialist: Behavioral support for a neural hypothesis. *American Naturalist* **151** 451–464.
[3] BEZUIDENHOUT, C. and GRIMMETT, G. (1991). Exponential decay for subcritical contact and percolation processes. *Ann. Probab.* **19** 984–1009. MR1112404
[4] BRAMSON, M. and DURRETT, R. (1988). A simple proof of the stability theorem of Gray and Griffeath. *Probab. Theory Related Fields* **80** 293–298. MR0968822
[5] BROWN, J. S. (1990). Habitat selection as an evolutionary game. *Evolution* **44** 732–746.
[6] CHESSON, P., PACALA, S. and NEUHAUSER, C. (2002). Environmental niches and ecosystem functioning. In *The Functional Consequences of Biodiversity* (A. P. Kinzig, S. W. Pacala, and D. Tilman, eds.) 213–245. Princeton Univ. Press.
[7] DURRETT, R. (1984). Oriented percolation in two dimensions. *Ann. Probab.* **12** 999–1040. MR0757768
[8] DURRETT, R. (1991). A new method for proving the existence of phase transitions. In *Spatial Stochastic Processes* 141–169. Birkhäuser, Boston. MR1144095
[9] DURRETT, R. and LEVIN, S. (1994). The importance of being discrete (and spatial). *Theor. Popul. Biol.* **46** 363–394.
[10] DURRETT, R. and NEUHAUSER, C. (1994). Particle systems and reaction-diffusion equations. *Ann. Probab.* **22** 289–333. MR1258879
[11] DURRETT, R. and NEUHAUSER, C. (1997). Coexistence results for some competition models. *Ann. Appl. Probab.* **7** 10–45. MR1428748
[12] ELTON, C. (1966). *Animal Ecology*. Sidgwick and Jackson, London.
[13] HARRIS, T. E. (1972). Nearest neighbor Markov interaction processes on multidimensional lattices. *Adv. Math.* **9** 66–89. MR0307392
[14] HARRIS, T. E. (1976). On a class of set valued Markov processes. *Ann. Probab.* **4** 175–194. MR0400468
[15] HUTCHINSON, G. E. (1957). Concluding remarks. *Cold Spring Harbour Symposium on Quantitative Biology* **22** 415–427.
[16] LIGGETT, T. M. (1985). *Interacting Particle Systems.* Springer, New York. MR0776231
[17] LIGGETT, T. (1999). *Stochastic Interacting Systems*: *Contact, Voter and Exclusion Processes.* Springer, Berlin. MR1717346
[18] NEUHAUSER, C. (1992). Ergodic theorems for the multitype contact process. *Probab. Theory Related Fields* **91** 467–506. MR1151806
[19] NEUHAUSER, C. and PACALA, S. (1999). An explicitly spatial version of the Lotka–Volterra model with interspecific competition. *Ann. Appl. Probab.* **9** 1226–1259. MR1728561
[20] NEUHAUSER, C. and LANCHIER, N. (2005). Coarseness of habitat affects evolution of specialization. Unpublished manuscript.

[21] Rosenzweig, M. L. (1987). Habitat selection as a source of biological diversity. *Evol. Ecol.* **1** 315–330.
[22] Woolhouse, M. E. J., Taylor, L. H. and Haydon, D. T. (2001). Population biology of multihost pathogens. *Science* **292** 1109–1112.

Department of Ecology,
 Evolution and Behavior
University of Minnesota
1987 Upper Buford Circle
St. Paul, Minnesota 55108
USA
E-mail: nicolas.lanchier@univ-rouen.fr
        neuha001@umn.edu